\documentclass[11pt]{amsart}

\newcommand{\newsection}[1]{\setcounter{equation}{0} \section{#1}}
\newtheorem{Theorem}{\sc Theorem}
\newtheorem{Lemma}[Theorem]{\sc Lemma}

\newtheorem{Corollary}[Theorem]{\sc Corollary}
\theoremstyle{definition}

\theoremstyle{remark}
\newtheorem{Remark}[Theorem]{\sc Remark}

\newcommand{\bea}{\begin{eqnarray}}
\newcommand{\eea}{\end{eqnarray}}
\newcommand{\be}{\begin{equation}}
\newcommand{\ee}{\end{equation}}
\newcommand{\ben}{\begin{eqnarray*}}
\newcommand{\een}{\end{eqnarray*}}
\def \qed {\hfill \vrule height6pt width 6pt depth 0pt}

\newcommand{\bal}{{\boldsymbol{\alpha}}}
\def\zbar{{\underline{z}}}
\def\wbar{{\underline{w}}}

\def\ul{\underline}%

\newcommand{\inner}[2]{\langle #1,#2 \rangle }

\newcommand{\scH}{\mathcal{H}}
\newcommand{\scM}{\mathcal{M}}

\newcommand{\scA}{\mathcal{A}}
\newcommand{\scK}{\mathcal{K}}

\def\C{\mbox{${\mathbb C}$}}
\def\A{\mbox{${\mathbb A}$}}
\def\D{\mbox{${\mathbb D}$}}
\def\T{\mbox{${\mathbb T}$}}

\newcommand{\cla}{\mathcal{A}}
\newcommand{\clb}{\mathcal{B}}

\newcommand{\clh}{\mathcal{H}}

\newcommand{\clm}{\mathcal{M}}

\newcommand{\bh}{\mathcal{B}(\mathcal{H})}

\begin{document}%
\title[Completely contractive maps]{Contractive and completely contractive
maps over planar algebras}
\author[Bhattacharyya]{Tirthankar Bhattacharyya}

\address{Department of Mathematics, Indian Institute of Science,
Bangalore 560012, India}

\email{tirtha@math.iisc.ernet.in}

\author[Misra]{Gadadhar Misra}
\address{Indian Statistical Institute, R. V. College Post, Bangalore
  560 059, India}

\email{gm@isibang.ac.in}

\thanks{The first named author acknowledges
  the support from Department of Science and Technology, India, Grant
  \# SR/ FTP/ MS-16/ 2001. The second named author acknowledges the
  support from the Indo-French Centre for the Promotion of Advanced
  Research, Grant \# IFC/2301-C/99/2396.}


\title[Contractive and completely contractive homomorphisms]
{Contractive and completely contractive homomorphisms of planar algebras}

\begin{abstract}
We consider contractive homomorphisms of a planar algebra ${\mathcal
A}(\Omega)$ over a finitely connected bounded domain $\Omega \subseteq
\C$ and ask if they are necessarily completely contractive.  We show
that a homomorphism $\rho:{\mathcal A}(\Omega) \to {\mathcal
B}(\mathcal H)$ for which $\dim({\mathcal A}(\Omega)/\ker \rho) = 2$
is the direct integral of homomorphisms $\rho_T$ induced by operators
on two dimensional Hilbert spaces via a suitable functional calculus
$\rho_T: f \mapsto f(T),~f\in {\mathcal A}(\Omega)$.  It is well-known that
contractive homomorphisms $\rho_T$, induced by a linear transformation
$T:\C^2 \to \C^2$ are necessarily completely contractive.
Consequently, using Arveson's dilation theorem for
completely contractive homomorphisms, one concludes that such a homomorphism $\rho_T$
possesses a dilation. In this paper, we construct this dilation explicitly.
In view of recent examples discovered by Dritschel and McCullough, we know that
not all contractive homomorphisms $\rho_T$ are completely contractive even if $T$
is a linear transformation on a finite-dimensional Hilbert space.
We show that one may be able to produce an example of a contractive homomorphism
$\rho_T$ of ${\mathcal A}(\Omega)$ which is not completely contractive
if an operator space which is naturally associated with the problem is not the MAX
space. Finally, within a certain special class of contractive
homomorphisms $\rho_T$ of the planar algebra ${\mathcal A}(\Omega)$,
we construct a dilation.
\end{abstract}
\maketitle

\newsection{Introduction}

All our Hilbert spaces are over complex numbers and are assumed to be
separable. Let $T \in \mathcal{B}(\mathcal{H})$, the algebra of bounded
operators on $\mathcal{H}$. The operator $T$ induces a homomorphism
$\rho_T: p \mapsto p(T)$, where $p$ is a polynomial.  Equip the
polynomial ring with the supremum norm on the unit disc, that is,
$\|p\| = \sup\{|p(z)|:z\in \D\}$. A well-known inequality due to
von Neumann (cf. \cite{von}) asserts that $\rho_T$ is contractive,
that is, $\|\rho_T\| \leq 1$ if and only if the operator $T$ is a
contraction. Thus in this case, contractivity of the homomorphism
$\rho_T$ is equivalent to the operator $T$ being a contraction.
As is well known, Sz.-Nagy \cite{nagy} showed that a contraction
$T$ on a Hilbert space $\mathcal{H}$ dilates to a unitary operator
$U$ on a Hilbert space $\mathcal{K}$ containing $\mathcal{H}$,
that is, $P p(U) h = p(T) h$ for all $h\in \mathcal{H}$ and any polynomial $p$,
where $P: \mathcal{K} \to \mathcal{H}$ is the projection of $\mathcal K$ onto
$\mathcal{H}$. The unitary
operator $U$ has a continuous functional calculus and hence induces a
$*$ - homomorphism $\varphi_U: C(\sigma(U)) \to \mathcal{B}(\mathcal{K})$.
It is easy to check that
$P[(\varphi_U)_{|\mathcal{A}(\mathbb{D})}(f)]_{|{\mathcal{H}}} =
\rho_T(f)$, for $f$ in ${\mathcal A}(\D)$, where ${\mathcal A}(\D)$
is the closure of the polynomials with respect to the supremum norm
on the disc $\D$.

Let $\Omega$ be a finitely connected bounded domain in $\C$.  We make
the standing assumption that the boundary of $\Omega$ is the disjoint
union of simple analytic closed curves. Let $T$ be a bounded linear operator on
the Hilbert space $\mathcal{H}$ with spectrum $\sigma(T) \subseteq
{\Omega}$.  Given a rational function $r = p/q $ with no poles in
the spectrum $\sigma(T)$, there is the natural functional calculus
$r(T) = p(T)q(T)^{-1}$. Thus $T$ induces a unital homomorphism $\rho_T
= r(T)$ on the algebra of rational functions ${\rm Rat}(\Omega)$ with
poles off $\Omega$.  Let $\mathcal{A}(\Omega)$ be the closure of ${\rm
Rat}(\Omega)$ with respect to the norm $\|r\| := \sup\{|r(z)\|:z\in
\Omega\}$.  Since functions holomorphic in a neighborhood of
$\bar{\Omega}$ can be approximated by rational functions with poles
off $\bar{\Omega}$, it follows that they belong to ${\mathcal A}
(\Omega)$.

The homomorphism $\rho_T$ is said to be {\em dilatable} if there
exists a normal operator $N$ on a Hilbert space $\mathcal{K} \supseteq
\mathcal{H}$ with $\sigma(N)\subseteq \partial\bar{\Omega}$ such
that the induced homomorphism $\varphi_N:C(\sigma(N)) \to
\mathcal{B}(\mathcal{K})$, via the functional calculus for the normal
operator $N$, satisfies the relation
\begin{equation} \label{dilation}
P( \varphi_N )_{|\mathcal{A}({\Omega})}(f) h = \rho_T(f) h,
\end{equation}
for $h$ in $\mathcal{H}$ and $f$ in $\mathcal{A}(\Omega)$.  Here
$P:\mathcal K \to \mathcal H$ is the projection of $\mathcal K$ onto
$\mathcal H$.

The observations about the disk prompt two basic questions:
\begin{enumerate}
\item[(i)] when is $\rho_T$ contractive;
\item[(ii)] do contractive homomorphisms $\rho_T$ necessarily dilate?
\end{enumerate}
For the disc algebra, the answer to the first question is given by
von Neumann's inequality while the answer to the second question
is affirmative -- Sz.-Nagy's dilation theorem. If the
domain $\Omega$ is simply connected these questions can be reduced
to that of the disc (cf. \cite{Dsar}).

If the domain $\Omega$ is the annulus, while no satisfactory
answer to the first question is known, the answer to the second
question was shown to be affirmative by Agler (cf. \cite{jim1}).

If $\rho_T:\mathcal{A}(\Omega) \to \mathcal{M}_2$ is a
homomorphism induced by an operator $T:\C^2 \to \C^2$ then it is
possible to obtain a characterization of contractivity and then
use it to show that the second question has an affirmative answer.
We do this in Section \ref{concrete}.  In Section \ref{rank2}, we
show that a larger class of contractive homomorphisms, we call
them contractive homomorphisms of {\em rank 2}, dilate. This is
done by proving that the rank $2$ homomorphisms are direct
integrals of homomorphisms induced by  two dimensional operators.

Arveson (cf. \cite{WBA2} and \cite{WBA3}) has shown that the existence
of a dilation of a contractive homomorphism $\rho$ of the algebra
$\scA(\Omega)$ is equivalent to complete contractivity of the
homomorphism $\rho$. We recall some of these notions in greater detail
in section \ref{opsp}.  We then show, how one may proceed to possibly
construct an example of a contractive homomorphism of the algebra
$\mathcal A(\Omega)$ which does not dilate.

In the final section of the paper, we obtain a general criterion for
contractivity.  This involves a factorization of a certain positive
definite kernel.  More importantly, we outline a scheme for
constructing the dilation of a homomorphism $\rho_T: \mathcal
A(\Omega) \to \mathcal M_n$ induced by an operator $T$ with distinct
eigenvalues.  This scheme is a generalization of the construction of
the dilation in section \ref{concrete}.

\newsection{Homomorphisms of Rank Two}\label{rank2}

A homomorphism $ \rho : \scA(\Omega) \rightarrow
\mathcal{B}(\mathcal{H})$ is said to be of {\em rank $n$} if it has
the property $\dim \big ( \scA(\Omega)/ \ker \rho \big )$ $=$ $n$.
In this section, we shall begin construction of dilation for
homomorphisms of rank 2. Nakazi and Takahashi showed that
contractive homomorphisms $\rho:\scA(\Omega) \to \mathcal{B}(\scH)$
of rank $2$ are completely contractive for any uniform sub-algebra
of the algebra of continuous functions $C(\bar{\Omega})$ (see
\cite{nt}). We would like to mention here that a generalization of this result 
was obtained by Meyer in Theorem 4.1 of \cite{ralf}. 
He showed that given a commutative unital closed subalgebra $\cla$ of $\mathcal B(\mathcal K)$
(for some Hilbert space $\mathcal K$) and a positive integer $d$, any $d-1$ contractive 
unital homomorphism $\rho : \cla \to \mathcal{M}_d$ is completely contractive.
In what follows, we construct explicit dilations for
homomorphisms from $ \scA(\Omega)$ to $ \mathcal{B}(\mathcal{H})$ of
rank two.

We first show that any homomorphism $\rho$ of rank $2$ is the direct
integral of homomorphisms of the form $\rho_T$ as defined in the
introduction, where $T \in \mathcal{M}_2$. The existence of dilation
of a contractive homomorphism $\rho_T$ induced by a two dimensional
operator $T$ is established in \cite{GM} by showing that the
homomorphism $\rho_T$ must be completely contractive.  It then
follows that every contractive homomorphism $\rho$ of rank $2$ must
be completely contractive. This implies by Arveson's theorem that
they possess a dilation. However, it is not always easy to construct
the dilation whose existence is guaranteed by the theorem of
Arveson. In this case, we shall explicitly construct the dilation of
a homomorphism of rank $2$. This is achieved by constructing the
dilation of a contractive homomorphism of the form $\rho_T$ for a
two dimensional operator $T$.

\begin{Lemma} \label{rank2lemma} If $\rho_T : \scA(\Omega) \rightarrow
\mathcal{B}(\mathcal{L})$ is a homomorphism of rank two, then up
to unitary equivalence, the Hilbert space $ \mathcal{L}$ is a
direct integral
$$ \mathcal{L} = \int^\oplus_\Lambda \mathcal{L}_\lambda d \nu(\lambda)\,, $$
where each $ \mathcal{L}_\lambda $ is two-dimensional. In this
decomposition, the operator $T$ is of the form
$$ T = \int^\oplus_\Lambda \left (
\begin{array}{cc}
z_1(\lambda) & c(\lambda) \\
0 & z_2(\lambda)
\end{array}\right ) d \nu(\lambda).$$
\end{Lemma}

\begin{proof}
To begin with, it is easy to see (see Lemma 1 of \cite{nt}) that
$\mathcal{L}$ is a direct sum of two Hilbert spaces ${\mathcal H}$
and ${\mathcal K}$ and the operator $T: \scH\oplus \scK \to
\scH\oplus \scK$ is of the form:
$$\left (\begin{array}{cc}
z_1 I_{\mathcal{H}} & C \\
0 & z_2I_{\mathcal{K}}
\end{array}\right ), \mbox{ with } z_1,\:z_2 \in \Omega \mbox{ or }
\left ( \begin{array}{cc}
z I_{\mathcal{H}} & C \\
0 & z I_{\mathcal{K}}
\end{array}\right ), \mbox{ with } z\in \Omega ,$$ where $C$ is a
bounded operator from $\mathcal{K}$ to $\mathcal{H}$. Now if we
put $\scK_0 = (\ker C)^\perp,~ \scK_1 = \ker C $, $\scH_0
=\overline{{\rm Ran}}\; C$ and $ ~\scH_1 = ({\rm Ran}\; C)^\perp$,
then with respect to the decomposition $\scK = \scK_0 \oplus
\scK_1$ and $\scH = \scH_0\oplus \scH_1$, we have
$$C = \left (\begin{array}{cc}
\tilde{C} & 0 \\
0 & 0
\end{array}\right ),$$
where the operator $\tilde{C}$ is from $\mathcal{K}_0$ to
$\mathcal{H}_0$. The polar decomposition of $\tilde{C}$ then
yields $\tilde{C} = V P$, where the operator $V$ is unitary and
$P$ is positive. We apply the spectral theorem to the positive
operator $P$ and conclude that there exists a unitary operator
$\Gamma:\int^\oplus_\Lambda \scH_\lambda d \nu(\lambda) \to \scK_0$
which intertwines the multiplication operator $M$ on the Hilbert
space $\int^\oplus_\Lambda \scH_\lambda d\nu(\lambda)$ and $P$.

Now notice that the operator $T:\scH\oplus \scK \to \scH\oplus
\scK$ can be rewritten as
$$
\begin{pmatrix}
z_1 I_{\mathcal{H}_1} & 0 & \tilde{C}_{\mathcal{K}_0 \to
\mathcal{H}_0}& 0 \\
0&z_1 I_{\mathcal{H}_0}&0&0\\
0&0&z_2I_{\mathcal{K}_0}&0\\
0 & 0& 0& z_2I_{\mathcal{K}_1}
\end{pmatrix}.
$$
Interchanging the third and the second column and then the second
and third row, which can be effected by a unitary operator, we see
that the operator $T$ is unitarily equivalent to the direct sum of
a diagonal operator $D$ and an operator $\tilde{T}$ of the form
$\left (
\begin{smallmatrix} z_1 I_{\mathcal{H}_0} & \tilde{C}_{\mathcal{K}_0
\to \mathcal{H}_0} \\ 0 & z_2I_{\mathcal{K}_0}
\end{smallmatrix}\right )$, where $\tilde{C}$ has dense range.  It is
clear that if we conjugate the operator $\tilde{T}$ by the operator
$I_{\mathcal{H}_0} \oplus U_{\mathcal{H}_0 \to \mathcal{K}_0}$, where
$U$ is any unitary operator identifying $\scH_0$ and $\scK_0$ then we
obtain a unitarily equivalent copy of $\tilde{T}$ (again, denoted by
$\tilde{T}$) which is of the form $\left ( \begin{smallmatrix} z_1
I_{\mathcal{H}_0} & \tilde{C}_{\mathcal{K}_0 \to \mathcal{H}_0}
U_{\mathcal{H}_0 \to \mathcal{K}_0} \\ 0 & z_2I_{\mathcal{K}}
\end{smallmatrix}\right )$.
Now, if we apply the polar decomposition to $\tilde{C}$ then we
see that the off diagonal entry is a positive operator on
$\scH_0$. One then sees that $\tilde{T}$ is unitarily equivalent
to $\left ( \begin{smallmatrix}
z_1 I_{\int^\oplus_\Lambda \scH_\lambda d\nu(\lambda)} & M \\
0 & z_2I_{\int^\oplus_\Lambda \scH_\lambda d\nu(\lambda)}
\end{smallmatrix}\right )$ via conjugation using the operator
$\Gamma\oplus \Gamma$.  We need to conjugate this operator one more
time using the unitary $W$ that identifies $\int^\oplus_\Lambda
\scH_\lambda d\nu(\lambda) \oplus \int^\oplus_\Lambda \scH_\lambda
d\nu(\lambda)$ and $\int^\oplus_\Lambda \scH_\lambda \oplus
\scH_\lambda d\nu(\lambda)$, where $W(s_1\oplus s_2)(\lambda) =
s_1(\lambda) \oplus s_2(\lambda)$ for $s_1\oplus_2 \in
\int^\oplus_\Lambda \scH_\lambda d\nu(\lambda)$. It is easy to
calculate $W \tilde{T} W^*$ and verify the claim.
\end{proof}

In view of the Lemma above, it is now enough to consider dilations
of homomorphisms $\rho_T$ where $T$ is a linear transformation on
$\C^2$.

\newsection{Dilations and Abrahamse-Nevanlinna-Pick Interpolation}
\subsection{} Consider any reproducing kernel Hilbert space
$\mathcal{H}_K$ of holomorphic functions on $\Omega$ with $K :
\Omega \times \Omega \to \C$ as the kernel. Assume that the
multiplication operator $M$ by the independent variable $z$ is
bounded. Then $M^* ( K( \cdot , z)) = \bar{z} K(\cdot ,z)$ and
it is clear by differentiation that
$M^* \bar{\partial}_z K( \cdot , z) = K(\cdot , z) +
\bar{z}\bar{\partial}_z K(\cdot ,z)$.

The matrix representation of the operator $M^*$ restricted to the
subspace $\mathcal M$ spanned by the two vectors $K(\cdot,z_1)$ and
$K(\cdot, z_2)$ has two distinct eigenvalues $\bar{z}_1$ and
$\bar{z}_2$. Similarly, the operator $M^*$ restricted to the
subspace $\mathcal N$ spanned by the two vectors $K(\cdot ,z)$ and
$\bar{\partial}_z K(\cdot , z)$ has only one eigenvalue $\bar{z}$ of
multiplicity $2$. In the lemma below, we identify certain $2$
dimensional subspaces of $\mathcal{H}_K \oplus \mathcal{H}_K$ which
are invariant under the multiplication operator $M^*$ and then
find out the form of the matrix. The reproducing kernel $K$
satisfies:

\begin{subequations} \label{reproduce}
\begin{eqnarray}
K(z_1,z_2) &=& \langle K(\cdot, z_2), K(\cdot,z_1)\rangle,\:z_1, z_2
\in \Omega,\\
(\partial_z K)( z , u) &=& \langle K(\cdot, u), \bar{\partial}_z K(
\cdot , z) \rangle, \:u,\, z\in \Omega. \label{differentiated}
\end{eqnarray} \end{subequations}
Using (\ref{reproduce}) and applying the Gram-Schmidt orthogonolization process
to the set $\{K(\cdot , z_1), K( \cdot , z_2)\}$, we get the orthonormal pair
of vectors
$$ e(z_1) = \frac{K ( \cdot , z_1)}{ K( z_1 , z_1)^{1/2}} \mbox{
and }   f (z_1 , z_{ 2})  =  \frac{ K ( z_1 , z_1) K(
\cdot , z_2) - K ( z_2 , z_2) K ( \cdot , z_1)}{ K ( z_1 , z_1)
^{1/2} \big ( K (z_1, z_1) K ( z_2 , z_2) - | K ( z_1 , z_2)
|^2\big )^{1/2}}.$$
Now for any $\mu\in\bar{\mathbb{D}}$, the pair of vectors
$$h_1(z_1 , z_2)  = \begin{pmatrix}
  0 \vspace*{2mm} \\  e (z_1)
\end{pmatrix} \mbox{ and }
h_2(z_1,z_2)  = \begin{pmatrix} (1 - | \mu |^2)^{1/2} e (z_2) \\
\mu f (z_1 , z_2)
\end{pmatrix}$$
are orthonormal in $\mathcal{H}_K \oplus \mathcal{H}_K$.
Similarly, using (\ref{differentiated}), orthonormalization of the
pair of vectors $\{K(\cdot , z), \bar{\partial}_z K( \cdot , z)\}$ produces the
orthonormal set $\{e(z), f(z)\}$, where
$$
e(z)  =  \frac{K( \cdot , z)}{ K( z , z)^{1/2}} \mbox{ and } 
f (z) =  \frac{  K ( z , z) \bar{\partial}_z K ( \cdot , z) -
\langle \bar{\partial}_z K ( \cdot , z) , K ( \cdot , z) \rangle K
( \cdot , z)}{ K ( z , z)^{1/2} \big ( K ( z , z) \,\|
\bar{\partial}_z K ( \cdot , z) \|^2  - | \langle \bar{\partial}_z
K ( \cdot , z) , K ( \cdot , z) \rangle |^2 \big )^{1/2} }.  
$$
and then for any $\lambda \in\bar{\mathbb{D}}$,
$$k_1(z)  = \begin{pmatrix}
  0 \vspace*{2mm} \\  e (z)
\end{pmatrix} \mbox{ and }
k_2(z)  = \begin{pmatrix} (1 - | \lambda |^2)^{1/2} e (z) \\
\lambda f (z)
\end{pmatrix}$$
form a set of two orthonormal vectors in $ \mathcal{H}_K \oplus
\mathcal{H}_K$.

Note that from the definition of $M^*$ it follows that $ M^* e
(z_1) = {\bar{z}_1} e (z_1) \mbox{ for all } z_1 \in \Omega.$
Therefore we have $(M^*\oplus M^*) h_1(z_1 , z_2) =
{\bar{z}_1} h_1(z_1 , z_2)$.  Now, \ben M^* f(z_1 , z_2) & = &
\frac{ K ( z_1 , z_1) {\bar{z}}_2 K_\alpha ( \cdot , z_2) - K (
z_2 , z_2) {\bar{z}}_1 K ( \cdot , z_1)}{ K ( z_1 , z_1) ^{1/2}
( K(z_1, z_1) K ( z_2 , z_2) - | K ( z_1 , z_2) |^2)^{1/2}} \\
& = & {\bar{z}}_2 f (z_1 , z_2) + \frac{({\bar{z}}_2 -
{\bar{z}}_1) K(z_1 , z_2)} {(K(z_1 , z_1) K(z_2 , z_2) - |K(z_1 ,
z_2)|^2)^{1/2}} e(z_1). \een  It follows that $\mathcal M$ is
invariant under $M^*\oplus M^*$. In particular, we have \ben
(M^*\oplus M^*) h_2(z_1,z_2) & = &
\begin{pmatrix} (1 - | \mu |^2)^{1/2} M^* e (z_2) \\
\mu M^* f (z_1 , z_2)
\end{pmatrix} \\
& = & \begin{pmatrix} (1 - | \mu |^2)^{1/2} {\bar{z}}_2 e (z_2) \\
\mu \big ( {\bar{z}}_2  f (z_1 , z_2) + \frac{({\bar{z}}_2 -
{\bar{z}}_1) K(z_1 , z_2)}
  {(K(z_1 , z_1) K(z_2 , z_2) - |K(z_1 , z_2)|^2)^{1/2}} e (z_1) \big )
\end{pmatrix} \\
&=& \bar{z}_2 \begin{pmatrix} (1 - | \mu |^2)^{1/2}  e (z_2) \\ \mu f (z_1 , z_2) \end{pmatrix}
+ \begin{pmatrix}0\\ \mu \frac{({\bar{z}}_2 -
{\bar{z}}_1) K(z_1 , z_2)}
  {(K(z_1 , z_1) K(z_2 , z_2) - |K(z_1 , z_2)|^2)^{1/2}} e (z_1) \end{pmatrix} \\
&=& \bar{z}_2 h_2(z_1,z_2) + \mu \frac{({\bar{z}}_2 - {\bar{z}}_1) |K(z_1 , z_2)|}
  {(K(z_1 , z_1) K(z_2 , z_2) - |K(z_1 , z_2)|^2)^{1/2}} h_1(z_1,z_2), 
\een where we have absorbed the argument of $K(z_1,z_2)$ in $\mu$.  

Now recall that $(M^* - \bar{z} ) K (\cdot ,z) = 0$.
Differentiating with respect to $ \bar{z}$, we obtain, $ M^*
\bar{\partial}_z K(\cdot , z)  = K(\cdot , z) + \bar{z}
\bar{\partial}_z K(\cdot , z).$ Thus the subspace $\mathcal N$
spanned by the vectors $k_1(z), k_2(z)$ is invariant under $M^*$.
A little more computation, similar to the one above,
gives us the matrix representation of the restriction of the operator
$M^*\oplus M^*$ to the subspace $\mathcal N$.  

So, we have proved the following Lemma.

\begin{Lemma} \label{subsp} The
two-dimensional space $\mathcal M$ spanned by the two vectors
$h_1(z_1 , z_2), h_2(z_1 , z_2)$ is an invariant subspace for the
operator $M^* \oplus M^*$ on $ \mathcal{H}_K \oplus
\mathcal{H}_K$ and the restriction of this operator to the
subspace $\mathcal M$ has the matrix representation
$$\begin{pmatrix}
{\bar{z}}_{1} & \frac{\mu ({\bar{z}}_2 - {\bar{z}}_1) |K(z_1 ,
z_2)|}{(K(z_1 , z_1) K(z_2 , z_2) - |K(z_1 , z_2)|^2)^{1/2}}
\vspace*{2mm} \\
0 & {\bar{z}}_{2}
\end{pmatrix}.$$
Similarly, the two-dimensional space $\mathcal N$ spanned by the
two vectors $k_1(z), k_2(z)$ is an invariant subspace for the
operator $M^* \oplus M^*$ on $ \mathcal{H} \oplus \mathcal{H}$
and the restriction of this operator to the subspace $\mathcal N$
has the matrix representation
$$\begin{pmatrix}
{\bar{z}} & \frac{\lambda K(z , z)} {(K(z , z) \| \bar{\partial}_z
K( \cdot , z) \|^2 - | \langle \bar{\partial}_z K ( \cdot , z),
K ( \cdot , z) \rangle |^2)^{1/2}} \vspace*{2mm} \\
0 & {\bar{z}}
\end{pmatrix}.$$
\end{Lemma}

Let $\mu, \lambda$ be a pair of complex numbers and fix a pair of
$2\times 2$ matrices $A_s$ and $B_t$ --
\begin{equation}
A_s = \begin{pmatrix} z_1 & 0
\\ s \mu (z_1-z_2)& z_2 \end{pmatrix},\: z_1,\,z_2\in \Omega
\mbox{ and } B_t = \begin{pmatrix} z & 0 \\ t \lambda & z
\end{pmatrix},\: z\in \Omega,
\end{equation}
where $s,t$ are a pair of positive real numbers. If we choose

\begin{subequations} \begin{eqnarray}
s:= s_K &=& \frac{|K(z_1 , z_2)|}{(K(z_1 , z_1) K(z_2 , z_2) - |K(z_1
, z_2)|^2)^{1/2}}, \mbox{ and }\\
t:= t_K &=& \frac{ K(z , z)} {(K(z , z) \| \bar{\partial}_z
K( \cdot , z) \|^2 - | \langle \bar{\partial}_z K ( \cdot , z),
K ( \cdot , z) \rangle |^2)^{1/2}},
\end{eqnarray} \end{subequations}
then it follows from the Lemma that the matrix $A_s$ (respectively, $B_t$) is the
compression of the operator $M \oplus M$ on the
Hilbert space $\scH_K\oplus\scH_K$ to the two dimensional
subspaces $\scM$  (respectively, ${\mathcal N}$) if and only if
$|\mu| \leq 1$ (respectively, $|\lambda| \leq 1$).

A natural family of Hilbert spaces $H^2_\bal(\Omega)$ consisting of
modulus automorphic holomorphic functions on $\Omega$ was studied in
the paper \cite{AD}.  This family is indexed by $\boldsymbol{\alpha}
\in \T^m$, where $m$ is the number of bounded connected components in
$\C \setminus \Omega$ and $\T$ is the unit circle. Each
$H^2_\bal(\Omega)$ has a reproducing kernel $K_\bal(z , w)$.  It was
shown in \cite{AD} that every pure subnormal operator with spectrum
$\bar{\Omega}$ and the spectrum of the normal extension contained in
$\partial \bar{\Omega}$ is unitarily equivalent to $M$ on one of
these Hilbert spaces.

In the following subsection, we will show that any contractive
homomorphism of the algebra ${\rm Rat}(\Omega)$ is of the form
$\rho_{A_s}$ or $\rho_{B_t}$ with $K=K_{\boldsymbol{\alpha}}$
and $|\mu| \leq 1$ and $|\lambda| \leq 1$ respectively. Since the
operator $M \oplus M$ is subnormal, we would have exhibited
the dilation.

\subsection{\sc Construction of Dilations}
\label{concrete}
The generalization of Nevanlinna-Pick theorem due to
Abrahamse states that given $n$ points $w_1, w_2, \ldots , w_n$ in
the open unit disk, there is a holomorphic function $f : \Omega
\rightarrow \C$ with $f(z_i) = w_i$ for $i=1,2, \ldots , n$ if and
only if the matrix
\begin{equation}
M(\ul{w}, \boldsymbol{\alpha}) \stackrel{\rm def}{=} \big (\!
\!\big ( (1-w_i\bar{w}_j) K_\bal(z_i, z_j) \big ) \!\!\big )
\end{equation}
is positive semidefinite. A deep result due to Widom (cf. \cite[page
140]{fisher}) shows that the map $\boldsymbol{\alpha} \mapsto
K_\bal(z ,w )$ is continuous for any fixed pair $(z , w)$ in $\Omega
\times \Omega$.

In what follows, we shall first show that a homomorphism $\rho:
\mbox{Rat}(\Omega) \to \scM_2$ is contractive if and only if it is of
the form $\rho_{A_s}$ or $\rho_{B_t}$ with $|\mu| \leq 1$ and $|\lambda| \leq 1$, respectively
and
\begin{subequations} \begin{equation}
s = s_\Omega(z_1,z_2) := \sup \{ |r(z_1)|^2 : r \in \mbox{ Rat}(\Omega),\: \|r\| \leq 1
\mbox{ and } r(z_2) = 0 \}
\end{equation}
for any fixed but arbitrary pair $z_1, z_2 \in \Omega$;
\begin{equation}
t = t_\Omega(z) := \sup \{ |r'(z) | : r \in
\mbox{ Rat}(\Omega),\:\|r\| \leq 1\mbox{ and } r(z) = 0\}
\end{equation}
\end{subequations}
for  $z\in \Omega$.

We wish to point out that the extremal quantities $s_\Omega(z_1,z_2)$ and
$t_\Omega(z)$ would remain the same even if we were to replace the
${\rm Rat}(\Omega)$ by the holomorphic function on $\Omega$.  The
solution to the first extremal problem, with holomorphic functions in
place of ${\rm Rat}(\Omega)$, exist by a normal family argument.
Let $F: \Omega\to \D$ be a holomorphic function
with $F(z_2) =0$ and $F(z_2) = a$, where we have set
$a = s_\Omega(z_1,z_2)$, temporarily.  It then follows that
$M((0,a), \boldsymbol{\alpha})$
must be non negative definite for all $\bal \in \T^m$.  Consequently, we have
$$
\det \begin{pmatrix}
K_\bal(z_1,z_1) & K_\bal(z_1,z_2) \\
K_\bal(z_2,z_1) & (1-a^2) K_\bal(z_2,z_2)
\end{pmatrix} \geq 0\\
$$
for all $\bal \in \T^m$.  This condition is
equivalent to requiring
\begin{equation}
|a|^2 \leq 1 - \frac{|K_\bal(z_1,z_2)|^2}{K_\bal(z_1,z_1) K_\bal(z_2,z_2)} \leq 1 -
\sup \big \{\frac{|K_\bal(z_1,z_2)|^2}{K_\bal(z_1,z_1) K_\bal(z_2,z_2)}: \bal \in \T^m \big \}.
\end{equation}
As we have pointed out earlier, since $\bal \to K_\bal(z_i,z_j)$ is continuous for any pair 
of fixed indices $i$ and $j$, there exists
a single $\bal_0$ depending only on $z_1, z_2$ for which the supremum in the above inequality
is attained. For this choice of $\bal_0$ and $a$, clearly the determinant of
$M((0,a), \boldsymbol{\alpha}_0)$ is zero.  It follows from \cite[Theorem 4.4, pp. 135]{fisher}
that the solution is unique and hence is a Blaschke product \cite[Theorem 4.1, pp. 130]{fisher}.

Similarly, the solution to the second extremal problem, with holomorphic
functions in place of ${\rm Rat}(\Omega)$, is a function which is holomorphic
in a neighborhood of $\bar{\Omega}$ \cite[Theorem 1.6, pp. 114]{fisher}.
Hence it is the limit of functions from ${\rm Rat}(\Omega)$. The following
Lemma first appeared as \cite[Remark 2, pp. 308]{GM}.

\begin{Lemma} \label{zero}
The homomorphism $\rho_{A_s}$  is contractive if and only if $\|r(A_s)\| \leq 1$
for all $r$ in ${\rm Rat}(\Omega)$ with $\|r\| \leq 1$ and $r(z_1) = 0$.

The homomorphism $\rho_{B_t}$  is contractive if and only if $\|r(B_t)\| \leq  1$
for all  $r$ in ${\rm Rat}(\Omega)$
with $\|r\| \le 1$ and $r(z) = 0$.
\end{Lemma}

\noindent{\sc Proof:} The two proofs are similar, so we shall
prove only (1). Suppose $ r( A)$ is a contraction for all $r \in \mbox{
Rat}(\Omega)$ with $\|r\| \leq 1$ and $r(z_1) = 0$. We have to
prove $ r( A)$ is a contraction for all $ r \in$ Rat$(\Omega)$
with $\|r\| \leq 1$.  For any such rational function $r$, let $r(z) = u$.
Put $\varphi_u (z) = \frac{z-u}{1 - \overline{u} z}$ and $\psi(z) = \varphi_u
( r(z))$. Then $\psi$ is in ${\rm Rat}(\Omega)$, $\|\psi\|\leq 1$ and
$\psi(z) = 0$. By hypothesis, $\psi(A)$ is a contraction. Now note that
$\varphi_u^{-1} (z) = \frac{z+u}{1 + \overline{u} z}$. Since
$\varphi_u^{-1}$ maps $ \mathbb{D}$ into $ \mathbb{D}$, by von
Neumann's inequality, $ \| r(A) \| = \| \varphi_u^{-1} \psi(A) \| \leq
1$. \qed

This lemma makes it somewhat simple to derive the contractivity conditions
for the homomorphisms induced by $A_s$ and $B_t$.

\begin{Lemma}
The homomorphism $\rho_{A_s}$  is contractive if and
only if  $s^2 = s_\Omega(z_1,z_2)^{-1} -1$ and $| \mu | \le 1$.
Similarly, the homomorphism $\rho_{B_t}$ is
contractive if and only if $t= t_\Omega(z)^{-1}$ and $| \lambda | \le 1$.
\end{Lemma}

\noindent{\sc Proof:} First, using the functional calculus for $A_s$, we see that
$$ r \begin{pmatrix}
  z_1 & 0 \\
  s \mu (z_1 - z_2) & z_{2}
\end{pmatrix}  = \begin{pmatrix}
  r(z_1) & 0 \\
 s \mu (r(z_1) - r(z_2)) & r(z_{2})
\end{pmatrix} = \begin{pmatrix}
  r(z_1) & 0 \\
  s \mu r(z_1) & 0
\end{pmatrix},$$
assuming $r(z_2) = 0.$
Therefore, contractivity of $\rho_{A_s}$ would imply
$$ s^2|\mu|^2 + 1 \leq \big (\sup\{|r(z_2)|^2 : r \in {\rm Rat} (\Omega),\,
\|r\| \leq 1{\rm~and~} r(z_2) = 0\} \big )^{-1} = s_\Omega(z_1,z_2)^{-1}.
$$
Or, equivalently, if we put $s = s_\Omega(z_1,z_2)^{-1} - 1$ then we must have $|\mu| \leq 1$.
Now an application of Lemma \ref{zero} completes the proof.

To obtain the contractivity condition for $\rho_{B_t}$, using the functional
calculus, we see that
$$ r \begin{pmatrix}
 z  & 0 \\
  t \lambda & z
\end{pmatrix} = \begin{pmatrix}
  r(z) &  0 \\
  t \lambda r'(z) & r(z)
\end{pmatrix} = \begin{pmatrix}
  0 & 0 \\
  a \lambda r'(z) & 0
  \end{pmatrix}$$
assuming $r(z)=0$.

Therefore, contractivity of $\rho_{B_t}$ would imply that
$$ t |\lambda| \le  (\sup \{ |r'(w) | : r \in \mbox{ Rat}( \Omega),\, \|r\| \leq 1
\mbox{ and } r(w) = 0\})^{-1}=t_\Omega(z)^{-1}.$$
Or equivalently, if we put $t=t_\Omega(z)^{-1}$ then we must have $|\lambda | \leq 1$.
\qed

We now have enough material to construct the dilation for a
homomorphism $ \rho_T : \mathcal{A}(\Omega) \rightarrow
\mathcal M_2$. In this case, $T$ is a $2 \times 2$ matrix with spectrum in
$\Omega$. Since we can apply a unitary conjugation to make $T$
upper-triangular, it is enough to exhibit the dilation for the two matrices
$T = A_s$ and $T=B_t$.

\subsection{\sc Dilation for $A_s$}
Recall that there exists an $\boldsymbol{\alpha}_0$ depending only on
$z_1$ and $z_2$ such that  $\det M((0,a), \boldsymbol{\alpha}_0) =0$. For now,
set $\bal_0 = \bal$. Let the subspace $\mathcal M$ of $ H^2_\bal
\oplus H^2_\bal$ be as in the first part of Lemma \ref{subsp}. For
brevity, let
$$m^2 =  1 - \frac{ | K_\bal (z_1 , z_2) |^2}{ K_\bal (z_1 , z_1)
K_\bal (z_2 , z_2) } > 0.$$  Then
$\det M((0,m) , \boldsymbol{\alpha}) = \begin{pmatrix}
  K_{\bal}(z_1, z_1)  & K_{\bal} (z_1 , z_2) \\
  K_{\bal} (z_2 , z_1) & (1 - m^2) K_{\bal} (z_2 ,z_2)
\end{pmatrix} = 0$
by definition of $m$.  As we have pointed out earlier, there is a holomorphic function
$f : \Omega \rightarrow \mathbb{D}$ such that $f(z_1) = 0$ and $f(z_2) = m$.
Moreover, if $g$ any holomorphic function from $\Omega$ to
$\mathbb{D}$ such that $g(z_1) = 0$, then the matrix $M((0,
g(z_2)), \boldsymbol{\alpha})$ is positive semidefinite, which
implies that 
$|g(z_2) |^2 \le 1 - \frac{ | K_\bal (z_1 , z_2) |^2}
{ K_\bal (z_1 , z_1)K_\bal (z_2 , z_2) }.$  Thus $m = \sup \{ |g(z_2)| : g
\mbox{ is a holomorphic function from } \Omega \mbox{ to }
\mathbb{D} \mbox{ and } g(z_1) = 0\}$. Hence
$$ s_\Omega(z_1,z_2)^2 = \frac{1}{m^2} - 1 =  \frac{| K_\bal(z_1,z_2)
|^2}{K_\bal(z_1,z_1) K_\bal(z_2,z_2) - | K_\bal(z_1,z_2) |^2}.
$$
So by the first part of Lemma \ref{subsp}, we have that the matrix of
restriction of the operator $M^*\oplus M^*$ to the subspace
$\mathcal M$ in the orthonormal basis $\{h_1(z_1 , z_2) , h_2(z_1,z_2)\}$
has the matrix representation $A_s^*$ with $s^2 = s_\Omega(z_1,z_2)^{-1} - 1$
whenever $|\mu| \leq 1$.

\subsection{} Having constructed the dilation, it is natural to find out
what the characteristic function is when $\Omega = \mathbb{D}$. In
this case, the general form of the matrix $T$ discussed above is \be
T :=
\begin{pmatrix}
  {z}_{1} & 0 \\
  {\mu} (1 - |z_{1}|^2)^{1/2} ( 1- |z_2|^2)^{1/2} & {z}_{2}
\end{pmatrix}. \label{Slambda} \ee
where $z_1$ and $z_2$ are two points in the open unit disk $
\mathbb{D}$ and $\mu \in \mathbb{C}$. We are using the
explicit value of $s_{\mathbb D}(z_1,z_2)$ for the unit disc.

\begin{Lemma} For $i=1,2$, let $\varphi_i (z) = (z - z_i)/(1 -
\overline{z}_i z)$. The characteristic function of $T$ is
$$ \theta_T (z) = \begin{pmatrix}
  (1 - |\mu|^2)^{1/2} \varphi_2(z) & -\mu \\
  \bar{\mu} \varphi_1(z) \varphi_2(z) & (1 - |\mu|^2)^{1/2}
  \varphi_1(z)
  \end{pmatrix}$$
\end{Lemma}

\noindent {\sc Proof:} Recall that $\mathcal M$ is the subspace
spanned by the orthonormal vectors $h_1(z_1 , z_2)$ and
$h_2(z_1,z_2)$. Since the compression of $M \oplus M$ to the
co-invariant subspace $\mathcal M$ is $T$, by Beurling-Lax-Halmos
theorem, we need to only find up to unitary coincidence (see
\cite{Nagy-Foias}, page 192 for definition) the inner function
whose range is $\mathcal{M}^\perp$. So let
$\begin{pmatrix} f \\
g
\end{pmatrix}$ be a vector in the orthogonal complement of $
\mathcal{M}$. The condition of orthogonality to $h_1$ implies that
$g(z_1) = 0$ which is equivalent to $ g = \varphi_1 \xi$ for
arbitrary $\xi \in H^2(\mathbb{D})$. Now the orthogonality
condition to $h_2$ implies that $ (1 - | \mu |^2)^{1/2} f(z_2)
+ \mu \xi (z_2) = 0$, which is the same as \be
\label{orthocond} (1 - | \mu |^2)^{1/2} \varphi_1(z_2) f(z_2)
+ \mu g(z_2) = 0.\ee This implies that there is an $\eta_1 \in
H^2( \mathbb{D} )$ such that
$$ (1 - | \mu |^2)^{1/2}  f + \mu g^\prime = \varphi_2 \eta_1.$$
It is obvious that conversely if $\begin{pmatrix} f \\ g
\end{pmatrix}$ is a function from $ H^2( \mathbb{D}) \oplus H^2(
\mathbb{D})$ such that $g$ is in range of $\varphi$ and satisfies
(\ref{orthocond}), then it is in the orthogonal complement of
$\mathcal{M}$.

Now let $\eta_2 = (1 - | \mu |^2)^{1/2} \xi - \bar{\mu}
f.$ Then
$$\theta \begin{pmatrix}
  \eta_{1} \\
  \eta_{2}
\end{pmatrix}  = \begin{pmatrix}
  (1 - |\mu |^2)^{1/2} \varphi_2 \eta_1 - \mu \eta_2 \\
  \bar{\mu} \varphi_1 \varphi_2 \eta_1 + (1 - |\mu |^2)^{1/2}
  \varphi_1 \eta_2
\end{pmatrix} = \begin{pmatrix}
  f \\ g
\end{pmatrix}.$$ Thus if $ \begin{pmatrix}
  f \\ g
\end{pmatrix}$ satisfies (\ref{orthocond}), then it is in the
range of $\theta$. Conversely, it is easy to see that any element
in the range of $\theta $ will satisfy (\ref{orthocond}). Thus the
orthogonal complement of $ \mathcal{M}$ in $\mathcal{H}$ is the
range of $\theta$. So $\theta$ is the characteristic function of
the given matrix. \qed

We would like to remark here that for $z_1 = z_2$, the
characteristic function $\theta_T(u)$ for $T :=
\Big ( \begin{smallmatrix}
  z & 0 \\ {\lambda} (1 - |z|^2) & z
\end{smallmatrix} \Big )$ can be obtained
directly from the definition in case $z=0$.
A little computation, using the transformation rule for
the characteristic function under a biholomorphic automorphism
of the unit disk \cite[pp. 239 - 240]{Nagy-Foias}, produces the
formula
$$
\theta_T(u) = \begin{pmatrix} (1-|\lambda|^2)^{1/2} \varphi(u) & \lambda
\cr \bar{\lambda} \varphi^2(u) & (1-|\lambda|^2)^{1/2} \varphi(u) \end{pmatrix},\:u\in \D
$$
in the general case.

Let $T_\mu$ be the matrix defined in (\ref{Slambda}).
Note that if $T_{\mu^\prime}$ and $T_\mu$ are two such matrices with
$|{\mu^\prime} | = |\mu |$, then \ben \theta_{T_{\mu^\prime}} (z) & = &
\begin{pmatrix}
  (1 - |{\mu^\prime}|^2)^{1/2} \varphi_2 & -{\mu^\prime} \\
  \bar{{\mu^\prime}} \varphi_1 \varphi_2 & (1 - |{\mu^\prime}|^2)^{1/2}
  \varphi_1
  \end{pmatrix} = \\
  & = & \begin{pmatrix}
  (1 - |\mu|^2)^{1/2} \varphi_2 & - e^{i\psi} \mu \\
  e^{-i\psi} \bar{\mu} \varphi_1 \varphi_2 & (1 - |\mu|^2)^{1/2}
  \varphi_1
  \end{pmatrix} \mbox{ for some } \psi \in [0,2\pi] \\
  & = & \begin{pmatrix} e^{i\psi/2} & 0 \\
  0 & e^{i\psi/2} \end{pmatrix} \begin{pmatrix}
  (1 - |\mu|^2)^{1/2} \varphi_2 & - \mu \\
   \bar{\mu} \varphi_1 \varphi_2 & (1 - |\mu|^2)^{1/2}
  \varphi_1
  \end{pmatrix} \begin{pmatrix} e^{i\psi/2} & 0 \\
  0 & e^{i\psi/2} \end{pmatrix} \\
  & = & \begin{pmatrix} e^{i\psi/2} & 0 \\
  0 & e^{i\psi/2} \end{pmatrix} \theta_{S_\mu} (z)
  \begin{pmatrix} e^{i\psi/2} & 0 \\
  0 & e^{i\psi/2} \end{pmatrix},\een
  and hence their characteristic functions coincide. So they are
  unitarily equivalent. Conversely, if $T_{\mu^\prime}$ and $T_\mu$ are
  unitarily equivalent, then their characteristic functions
  coincide and hence the singular values of the characteristic
  functions are same. Note that when $z_1 \neq z_2$, we have
  $$\theta_{T_{\mu^\prime}} (z_1) \theta_{T_{\mu^\prime}}
  (z_1)^* = \begin{pmatrix}
  (1 - |{\mu^\prime}|^2) |\omega|^2 + |{\mu^\prime} |^2 & 0 \\
  0 & 0 \end{pmatrix} $$
  for some $\omega \in \mathbb{C}$ (independent of ${\mu^\prime}$).
  When $z_1=z_2$, then
  $$\theta_{T_{\mu^\prime}} (z_1) \theta_{T_{\mu^\prime}}
  (z_1)^* = \begin{pmatrix}
  0 & |{\mu^\prime} |^2 \\
  0 & 0 \end{pmatrix} .$$ In either case, coincidence of
  $\theta_{T_{\mu^\prime}}$ and $\theta_{T_\mu}$ mean that $|{\mu^\prime} |
  = |\mu |$. Thus using the explicit characteristic function we
  have proved the following.

\begin{Theorem}
Two matrices $T_{\mu^\prime}$ and $T_\mu$ as defined in (\ref{Slambda})
are unitarily equivalent if and only if $|{\mu^\prime} | = |\mu |$.
\end{Theorem}

\subsection{\sc Dilation for $B_t$}
We now shift our attention to the construction of dilation when the
homomorphism $\rho_T$ is induced by a $2 \times 2$ matrix $T$ with
equal eigenvalues. So $\sigma(T) = \{z\}$. The domain $\Omega$ has its
associated Szego kernel which is denoted by $\hat{K}_\Omega(z,w)$. Recall that
a generalization due to Ahlfors to multiply connected domains of the Schwarz 
lemma says that
$$t_\Omega(z) := \big (\sup \{ |r'(z) | : r \in \mbox{Rat}(\Omega), \|r\| \leq 1 \mbox{ and } r(z)
= 0\}\big )^{-1} = \hat{K}_\Omega(z,z)^{-1}.$$   
Let $ \partial \Omega$ be the topological boundary of $\Omega$ and let
$|d\nu|$ be the arc-length measure on $\partial \Omega$. Consider
the measure $dm = | \hat{K}_\Omega(\nu,z)|^2 |d\nu|$, and let the
associated Hardy space $H^2(\Omega, dm)$ be denoted by $\mathcal{H}$.
The measure $dm$ is mutually absolutely continuous with respect to the
arc length measure.  Thus the evaluation functionals on $\mathcal{H}$
are bounded and hence $ \mathcal{H}$ possesses a reproducing kernel $K$.
Then it is known that $K$ satisfies the property:
$$ \frac{K(z , z)} {(K(z , z) \| \partial_{\bar{z}}
K( \cdot , z) \|^2 - | \langle \partial_{\bar{z}} K ( \cdot , z) , K
( \cdot , z) \rangle |^2)^{1/2}} =  
\hat{K}_\Omega(z,z)^{-1},$$ see \cite[Theorem 2.2]{GM}. Now a (subnormal)
dilation for $A_s := \begin{pmatrix}
  z &  0 \\
  \lambda s_\Omega(z) & z
\end{pmatrix}$, where $|\lambda|\leq 1$,
is the operator $ M \oplus  M$ on the Hilbert space $\mathcal{H}
\oplus \mathcal{H}$. This is easily verified since the restriction
of $M^*\oplus M^*$ to the subspace $\mathcal N$ which was
described in the second part of Lemma \ref{subsp} is $A_s^*$.

\begin{Remark}
If we choose $|\mu| = 1$ then the $A_s^*$ is the restriction
of $M^*$ to the two dimensional subspace spanned by the vectors
$K_\bal(\cdot, z_1)$ and $K_\bal(\cdot, z_2)$ in the Hardy space
$H^2_\bal(\Omega)$ by our construction.
Except in this case, the dilation of the homomorphism $\rho_{A_s}$ we have
constructed is a minimal subnormal dilation. (This dilation then may be extended to
a minimal normal dilation.)  While it is known that a minimal dilation is not
unique when $\Omega$ is finitely connected, our construction gives a measure of
this non-uniqueness. More explicitly, for each $\bal_0 \in \T^m$ for which
$$\sup \{\frac{|K_\bal(z_1,z_2)|^2}{K_\bal(z_1,z_1) K_\bal(z_2,z_2)}:
\bal \in \T^m \}=  \frac{|K_{\bal_0}(z_1,z_2)|^2}{K_{\bal_0}(z_1,z_1)
K_{\bal_0}(z_2,z_2)},
$$
the matrix representation of the operator $M^*\oplus M^*$ restricted to
the $2$ dimensional subspace $\scM$ of the Hilbert space
$H^2_{\bal_0} \oplus H^2_{\bal_0}$ equals $A_s$.
\end{Remark}

\newsection{ The Operator Space} \label{opsp}

The problem that we are considering naturally gives rise to an
operator space structure. In this section, we show that.  We begin
by recalling basic definitions.

A vector space $X$ is called an operator space if for each $k\in
\mathbb{N}$, there are norms $\|\cdot\|_k$ on $X \otimes
\mathcal{M}_k$ such that
\begin{enumerate}
\item whenever $A = (( a_{ij} )) \in \mathcal{M}_k, (\!(x_{ij})\!) \in
X \otimes \mathcal{M}_k$ and $ B = (( b_{ij} )) \in \mathcal{M}_k$,
then $$\|A \cdot (\!( x_{ij} )\!) \cdot B\|_k \leq \|A\|
\|(\!(x_{ij} )\!)\|_k \|B\|$$ where $A \cdot (\!( x_{ij} )\!) \cdot
B = (( \sum_{p=1}^m \sum_{l=1}^k a_{ip} x_{pl} b_{lj} )) \in X
\otimes \clm_{k}$ and $\|A\|$ and $\|B\|$ are operator norms on
$\clm_k = \clb( \mathbb{C}^k)$.
\item For all positive integers $m,k$ and for all $R \in X
\otimes \clm_k$ and $S \in X \otimes \clm_m$, we have
$$ \Big \| \Big ( \begin{smallmatrix} R & 0 \\
0 & S
\end{smallmatrix} \Big ) \Big \|_{m+k} = \max\{\|R\|_m, \|S\|_k\}.$$
\end{enumerate}

Two such operator spaces $(X, \| \cdot \|_{X,k})$ and $(Y, \| \cdot
\|_{Y,k})$ are said to be completely isometric if there is a linear
bijection $ \tau : X \to Y$ such that $\tau \otimes I_k: (X,  \|
\cdot\|_{X, k} ) \to (Y,  \| \cdot\|_{Y, k} )$ is an isometry for
every $k \in \mathbb{N}$.

Let $X$ be an operator space and let $\rho: X \to \clb(\mathcal{H})$
be a linear map, where $\clh$ is a Hilbert space. If for each $k\in
\mathbb{N}$, the map $\rho \otimes I_k: (X, \| \cdot \|_k) \to
\mathcal{B}(\mathcal{H}\otimes\clm_k)$ is contractive then $\rho$ is
said to be {\em completely contractive}. Let $\clh$ be
finite-dimensional, let $T \in \clb(\clh)$, let $X = \cla(\Omega)$
and let $\rho = \rho_T$ be as defined earlier. We assume that the
eigenvalues $z_1, z_2, \ldots, z_n$ of $T$ are distinct.

To begin with, we introduce a notation. We denote by $I_\zbar^k$ the
subset of $\mathbb C^n \otimes \mathcal M_k$ defined as  $$
I_\zbar^k  = \{ (R(z_1), R(z_2), \ldots , R(z_n)) : R \in \cla(
\Omega ) \otimes \mathcal M_k \mbox{ and } \|R\| \leq 1 \}$$ where
$\|R\| = \sup_{ z\in \overline{\Omega}} \|R(z)\|$. When $k=1$, we
denote it by $I_{\ul{z}}$ rather than $I_\zbar^1$.

\begin{Lemma}
The set $I_\zbar$ defined above is a compact set.
\end{Lemma}

\proof Clearly, $I_\zbar$ is a subset of $\bar{D}^n$. So it is
enough to show that $I_\zbar$ is a closed set. Recall from Section 3
that the generalization of Nevanlinna-Pick theorem due to Abrahamse
states that given $n$ points $w_1, w_2, \ldots , w_n$ in the open
unit disk, there is a holomorphic function $f : \Omega \rightarrow
\mathbb{C}$ with $f(z_i) = w_i$ for $i=1,2, \ldots , n$ if and only
if the matrix

\begin{equation}
M(\ul{w}, \boldsymbol{\alpha}) \stackrel{\rm def}{=} \big (\! \!\big
( (1-w_i\bar{w}_j) K_\alpha(z_i, z_j) \big ) \!\!\big )
\end{equation}
is positive semidefinite for all $\bal \in \mathbb{T}^m$. So

$\begin{array}{ccl}
  I_\zbar  & = & \{ (w_1, w_2, \ldots , w_n) \in \overline{\mathbb{D}}^n\,:\,
\mbox{ the matrix } M(\wbar , \boldsymbol{\alpha}) \mbox{ is
positive semidefinite for all } \boldsymbol{\alpha}\in \mathbb{T}^m \}\\
   & = & \{ (w_1, w_2, \ldots , w_n) \in \overline{\mathbb{D}}^n\,:\,
\lambda_{\rm min} ( M(\wbar , \boldsymbol{\alpha}) ) \ge 0 \mbox{
for all }
\boldsymbol{\alpha}\in \mathbb{T}^n \}\\
   & = & \cap_{ \bal \in \mathbb{T}^m} \{ (w_1, w_2, \ldots , w_n) \in
\bar{\mathbb{D}}^n\,:\, \lambda_{\rm min} ( M(\wbar ,
\boldsymbol{\alpha}) )
\ge 0 \} \\
   & = & \cap_{ \bal \in \mathbb{T}^m} \left( \lambda_{\rm min}
( M(\wbar , \boldsymbol{\alpha}) ) \right)^{-1} \left( [0, \infty )
\right)
\end{array}$
where $\lambda_{\rm min}(A)$ for a hermitian matrix $A$ denotes its
smallest eigenvalue. It is a continuous function on the set of
hermitian matrices (see for example, \cite[ Corollary
III.2.6]{bhatia}). Thus $ \wbar \to \lambda_{\rm min} ( M(\wbar ,
\boldsymbol{\alpha}) )$ is a continuous function on $\mathbb{C}^n$.
Since arbitrary intersection of closed sets is closed, $I_\zbar$ is
a closed set. \qed

It is easy to see that the set $I_\zbar^k$ is convex and balanced,
so it is the closed unit ball of some norm on $ \mathbb{C}^n \otimes
\mathcal M_k$. The sets of the form $I_\zbar^k$ were first studied,
in the case $k=1$, by Cole and Wermer \cite{CW}.
The sets $I_\zbar^k$ are examples of {\em matrcially hyperconvex}
sets studied by Paulsen in \cite{paul2}. Paulsen points out that the
sequence of sets $I_\zbar^k\subseteq \mathbb C^n \otimes \mathcal
M_k$ determines an operator space structure on $\mathbb C^m$, that
is, the set $I_\zbar^k$ determines a norm $\|\cdot\|_{\ul{z},k}$ in
$\mathbb{C}^n \otimes \mathcal M_k$ such that $I_\zbar^k$ is the
closed unit ball in this norm and the sequence $\{ \mathbb{C}^n
\otimes \mathcal M_k , \|\cdot\|_{\ul{z},k} \}$ satisfies the
conditions (1) and (2) above. We denote this operator space by
$\mathrm{HC}_{\Omega,\ul{z}}(\mathbb C^n)$. Paulsen also notes that
this operator space is completely isometric to a quotient of a
function algebra. Indeed, it is not difficult to see that
$\mathrm{HC}_{\Omega,\ul{z}}(\mathbb C^n)$ is completely
isometrically isomorphic to the quotient $C^*$-algebra
$C(\overline{\Omega}) / \mathcal{Z}$ where $\mathcal{Z} = \{ f \in
\cla(\Omega) : f(z_1) = f(z_2) = \cdots = f(z_n) = 0\}$.
If $k=1$, we will write $\|\cdot\|_\zbar$ rather than
$\|\cdot\|_{\zbar,1}$.

\begin{Lemma}
There are $n$ matrices $V_1, V_2, \ldots ,V_n \in \clm_n$  such that
the map $\rho_T \otimes I_k : \mathcal{A}(\Omega) \otimes \clm_k
\rightarrow \clm_n \otimes \clm_k $ is of the form
$$ (\rho_T \otimes I_k) R = V_1 \otimes R(z_1) + V_2 \otimes R(z_2) +
\cdots + V_n \otimes R(z_n)$$ for any $R \in \mathcal{A}(\Omega)
\otimes \clm_k$ and any $k \in \mathbb{N}$. The matrices $V_i$
depend on the set $\{z_1, z_2, \ldots , z_n\}$.
\end{Lemma}

\noindent {\sc Proof}: If $F$ and $G$ are two elements of
$\mathcal{A}(\Omega) \otimes \clm_k$ which agree on the set $\{z_1,
z_2, \ldots ,z_n\}$, then define $H \in \mathcal{A}(\Omega) \otimes
\clm_k$ by $H = F - G$. Then $H$ vanishes at the points $z_1, z_2,
\ldots z_n$ and hence $H(z) =  (z - z_1)(z - z_2)\ldots (z - z_n)
W(z)$ for some $W$ in $\mathcal{A}(\Omega) \otimes \clm_k$.  By the
functional calculus,
$$(\rho_T \otimes I_k) H = (T - z_1) (T - z_2) \ldots (T -
z_n) W(T).$$ Note that $(z - z_1) (z - z_2) \ldots (z - z_n)$ is the
characteristic polynomial of $T$ and by Cayley-Hamilton theorem, $
(T - z_1)(T - z_2)\ldots (T - z_n) = 0$. Thus $(\rho_T \otimes I_k)H
= 0$. So if $F, G \in \mathcal{A}(\Omega) \otimes \clm_k$ are such
that $F(z_i) = G(z_i)$ for all $i=1,2, \ldots , n$, then $(\rho_T
\otimes I_k)F = (\rho_T \otimes I_k) G$.
Now define $V_1, V_2, \ldots ,V_n$ by
$$ V_i = \rho_T \left( \frac{(z - z_1) \ldots (z - z_{i-1})(z - z_{i+1}) \ldots
 (z - z_n)}{(z_i- z_1) \ldots (z_i - z_{i-1})(z_i - z_{i+1}) \ldots
 (z_i - z_n)}\right)$$ for $i=1,2, \ldots ,n$. Given $R \in
\mathcal{A}(\Omega) \otimes \clm_k$, it agrees with the function
$$\tilde{R}(z) = \sum_{i=1}^n\frac{(z - z_1) \ldots (z - z_{i-1})(z - z_{i+1})\ldots
 (z - z_n)}{(z_i- z_1) \ldots (z_i - z_{i-1})(z_i - z_{i+1}) \ldots
 (z_i - z_n)} R(z_i)$$ on the set $\{z_1, z_2, \ldots ,z_n\}$ and
hence 
\ben (\rho_T \otimes I_k)R  &=& (\rho_T \otimes I_k)\tilde{R} \\ 
& = & \sum_{i=1}^n \rho_T \left(\frac{(z - z_1) \ldots (z - z_{i-1})(z -
z_{i+1}) \ldots (z - z_n)}{(z_i- z_1)
\ldots (z_i - z_{i-1})(z_i - z_{i+1}) \ldots (z_i - z_n)}\right) \otimes R(z_i)\\
& = & \sum_{i=1}^n V_i \otimes R(z_i) \een 
completing the proof of the Lemma.  \qed 

At this point, we note that $\mathcal{A}(\Omega)$ being a 
closed sub-algebra of the commutative $C^*$-algebra of all
continuous functions on the boundary of $\Omega$ inherits a natural
operator space structure, denoted by ${\rm MIN}(\mathcal
A(\Omega))$. Recall that a celebrated theorem of Arveson says that a
contractive homomorphism $\rho_T :\mathcal{A} (\Omega) \to
\mathcal{L}(\mathcal{H)}$ dilates if and only if it is completely
contractive when $\mathcal A(\Omega)$ is equipped with the MIN
operator space structure. The contractivity and complete
contractivity of the homomorphism $\rho_T$ amount to respectively
\begin{equation} \label{lincont}
\sup \{\|w_1 V_1 + w_2 V_2 + \cdots + w_n V_n \| : \wbar =(w_1, w_2,
\ldots , w_n) \in I_\zbar  \} \leq 1
\end{equation}
where $\| \cdot \|$ is the operator norm on $ \clm_n$ and
\begin{equation}
\label{lincompcont} \sup \{\| \sum_{i=1}^n V_i \otimes W_i  \| : W_i
\in \clm_k \mbox{ and } W=(W_1, W_2, \ldots ,W_n) \in I_\zbar^k
\mbox{ for } k \ge 1 \} \leq 1
\end{equation}
where $\| \cdot \|$ is the operator norm on $ \clm_n \otimes
\clm_k$. Now, we state the following theorem whose proof is evident
from the discussion above.
\begin{Theorem} \label{homeqlin}
The contractive homomorphism $\rho_T: \mathcal{A}(\Omega) \to
\mathcal{M}_n$ is completely contractive with respect to the MIN
operator space structure on $\mathcal A(\Omega)$ if and only if the
contractive linear map $L_T:(\mathbb C^n,\|\cdot\|_\zbar) \to
\mathcal{M}_n$ defined by $L_T(\wbar) = w_1V_1 + w_2 V_2 + \cdots +
w_nV_n$ is completely contractive on the operator space ${\mathrm
{HC}}_{\Omega,\underline{z}}(\mathbb C^n)$.
\end{Theorem}

The theorem above brings us to our concluding remarks of this
section. Given a Banach space, there are natural operator space
structures on it, denoted by MAX($X$) and MIN($X$). We refer the reader to
\cite{vern} for definitions and basic details. However, this theorem
shows that if ${\mathrm {HC}}_{\Omega,\underline{z}}(\mathbb
C^n)$ was completely isometric to MAX($\mathbb C^n,
\|\cdot\|_\zbar$), then every contractive homomorphism $\rho_T$ of
the algebra $\mathcal A(\Omega)$, induced by an $n$ - dimensional
linear transformation $T$ with distinct eigenvalues in $\Omega$,
will necessarily dilate. This gives rise to the question of
determining when ${\mathrm {HC}}_{\Omega,\underline{z}}(\mathbb
C^n)$ is the same as MAX($\mathbb C^n ,\|\cdot\|_\zbar$) which is interesting in its
own right. Agler \cite{jim1} proved that all contractive
homomorphisms of the algebra $\mathcal{A}(\A)$, where $\A =\{z\in
\mathbb{C}: r< |z| < 1 \}\subseteq \mathbb{C}$ is the annulus for a
fixed  $r$ in $(0,1)$, are completely contractive. This implies that
$\mathrm{HC} _{\mathbb A,\underline{z}}(\mathbb{C}^n)$ is completely
isometric to MAX($\mathbb C^n,\|\cdot\|_\zbar$) for every $n \in \mathbb{N}$.

In \cite{vern}, Paulsen related a problem similar to the one that we
are considering to certain questions in the setting of operator spaces 
and thereby solved it.
For $n \ge 1$, let $G$ be a closed unit ball in $\mathbb{C}^n$ corresponding
to a norm $\| \cdot \|_G$ on $\mathbb{C}^n$. Let $\cla(G)$ denote the closure
of polynomials in $C(G)$, the algebra of all continuous functions on $G$
equipped with the sup norm. It is easy to see that there is a unital contractive
homomorphism $\rho : \cla(G) \to \bh$, for some Hilbert space $\clh$
which is not completely contractive if and only if
MIN$( \mathbb{C}^n , \| \cdot \|_G)$ is not completely isometric
to MAX$(\mathbb{C}^n , \| \cdot \|_G)$. Paulsen proved the remarkable result
that for $n \ge 5$,
\be \label{MINMAX} \mbox{MIN} ( \mathbb{C}^n , \| \cdot \|_G) \mbox{ is not
completely isometric to MAX} (\mathbb{C}^n , \| \cdot \|_G), \ee
for any closed unit ball $G$. For $n=2$, Ando's theorem implies that
MIN$( \mathbb{C}^2 , \| \cdot \|_{ \mathbb{D}^2})$ is completely isometric
to MAX$(\mathbb{C}^2 , \| \cdot \|_{ \mathbb{D}^2})$.  The fact that      
(\ref{MINMAX}) holds for $n \ge 3$ and any closed unit ball
$G$ is pointed out in \cite[Exercise 3.7]{Pisier}. 
In the same spirit, a similar question about a class of homomorphisms, 
first introduced by Parrott \cite{skp} (see also \cite{M}, \cite{MP} and \cite{MS}),   
led Paulsen to define a natural operator space which he called COT. 
Let $G$ be a unit ball and let $\wbar$ be a
point in the interior of $G$. Let $X$ be the Banach space
$X = ( \mathbb{C}^n , \| \cdot \|_{G,\wbar})$, where 
$\|\cdot\|_{G, \wbar}$ is the Caratheodory norm of $G$ at the point $\wbar$. 
The question of whether COT$_\wbar(X)$ is completely 
isometric to MIN($X^*$) for $\wbar \in G$ was first raised in
\cite{vern}.  He showed that the answer is affirmative when $\wbar =
0$. Later in an unpublished note, it was shown by Dash
\cite{mihir} that COT$_\wbar (G)$ and MIN$(X^*)$ are not 
necessarily completely isometric.  
The question of deciding whether a contractive homomorphism 
$\rho_T : \cla(\Omega) \to \bh$ 
is completeley contractive or not is similar in nature.  
It amounts to deciding if ${\mathrm {HC}}_{\Omega,\underline{z}}(\mathbb C^n)$
is completely isometric to MIN$(\mathbb{C}^n , \| \cdot \|_z)$ or 
MAX$(\mathbb{C}^n , \| \cdot \|_z)$. 
It is likely that the operator space ${\mathrm
{HC}}_{\Omega,\underline{z}}(\mathbb C^n)$ is completely isometric to
MIN$(\mathbb{C}^n , \| \cdot \|_z)$ for every $n \ge 3$. We pose this
as an open problem whose solution defies us at the moment.

\newsection{A Factorization Condition}
Let $T$ be a linear transformation on an $n$ dimensional Hilbert space 
space $V$ with distinct eigenvalues $z_1, z_2, \ldots , z_n$ in $\Omega$.
Let  $v_1, v_2, \ldots ,v_n$ be the $n$ linearly independent eigenvectors
of $T^*$. If $\sigma = \{ z_1, z_2, \ldots ,z_n\}$, then
define a positive definite function $K:\sigma \times \sigma \to
\C$ by setting 
\be \label{k} \big (\!\! \big ( K(z_j, z_i) \big )\!\! \big )_{i,j}^n 
:= \big (\!\!\big (\inner{v_i}{v_j} \big )\!\! \big )_{i,j=1}^n.
\ee 
As before, let $\rho_T: \mathcal A(\Omega) \to \mathcal L(V)$
be the homomorphism induced by $T$. Suppose there exists a
dilation of the homomorphism $\rho_T$. Then it follows from
\cite[Theorem 2]{AD1} that there is a flat unitary vector bundle
$\mathcal{E}$ of rank $n$ (see \cite{AD} for definitions and
complete results on model theory in multiply connected domains)
such that $\rho_T(f)$ is the compression of the subnormal operator $M_f$ 
on the Hardy space $H_\mathcal{E}^2(\Omega)$  
to a semi-invariant subspace in it. Consequently, there exists a homomorphism 
\be \label{generaldilation} {\rho}_M: \mathcal A(\Omega) 
 \to \mathcal{B}( H_\mathcal{E}^2(\Omega)) \ee 
dilating $\rho_T$.  The homomorphism ${\rho}_M$ is induced by the 
multiplication operator $M$ on $H_\mathcal{E}^2(\Omega)$ which is subnormal.
Thus the homomorphism $\rho_N: C(\partial \Omega) \to \mathcal B(\mathcal H)$ induced by
the normal extension $N$ on the Hilbert space $\mathcal H \supseteq H_\mathcal{E}^2(\Omega)$ 
of the operator $M$ is a dilation of the homomorphism $\rho_T$ in the sense of (\ref{dilation}). 
The multiplication operator $M$ on $H^2_\mathcal E$ is called a bundle shift. 
We recall \cite[Theorem 3]{AD} that $\dim \ker (M - z)^* = n$.  
Let $K^\mathcal E_z(i)$, $i=1,2, \ldots, n$ be a basis (not necessarily orthogonal) 
of $\ker (M - z)^*$.  We set 
\begin{equation} \label{hardyker}
K^\mathcal E(z_j, z_i) := 
\big ( \!\! \big ( \inner{K^\mathcal E_{z_i}(\ell)}{K^\mathcal E_{z_j}(p)}
\big )\!\! \big )_{\ell,p=1}^n,\: {\rm for~} 1\leq i,j \leq n.
\end{equation}

If $\rho_T$ dilates then the linear transformation $T$ can be realized
as the compression of the operator $M$ on $H_\mathcal{E}^2(\Omega)$ to
an $n$-dimensional co-invariant subspace, say $\mathfrak M \subseteq
H_\mathcal{E}^2(\Omega)$.  The subspace $\mathfrak M$ must consist of
eigenvectors of the bundle shift $M$.  Let $x_i$, $1\leq i \leq n$, be a
set of $n$ vectors in $\C^n$ and $\mathfrak M = \{\sum_{\ell=1}^n
x_i(\ell) K^\mathcal{E}_{z_i}(\ell): 1 \leq i \leq n\}$. The map which
sends $v_i$ to $\sum_{\ell=1}^n x_i(\ell) K^\mathcal{E}_{z_i}(\ell)$,
$1\leq i \leq n$, intertwines $T^*$ and the restriction of $M^*$ to 
$\mathfrak M$.  For this map to be an isometry as well, we must have
\begin{equation} \label{vtoE} 
\langle v_i, v_j \rangle 
= \langle K^\mathcal{E} (z_j , z_i) x_i , x_j \rangle,\: x_i\in \C^n,\:1\leq i\leq n.
\end{equation}
Conversely, if there is a flat unitary vector bundle
$\mathcal{E}$ and $n$ vectors $x_1, x_2, \ldots ,x_n$ in
$\C^n$ satisfying (\ref{vtoE}), then $\rho_T$ obviously
dilates. So we have proved the following theorem.

\begin{Theorem} The homomorphism $\rho_T$ is dilatable to a homomorphism
$\tilde{\rho}$ if and only if the kernel $K$, as defined in
(\ref{k}), can be written as  
$$ K(z_j , z_i) = \langle K^\mathcal{E} (z_j , z_i) x_i , x_j \rangle, 
\:{\rm ~for~ some~ choice~ of~}x_1, \ldots ,x_n\in \C^n,$$ 
where $K^\mathcal{E}(z_i,z_j)$ is defined in (\ref{hardyker}). 
\end{Theorem}

It is interesting to see how contractivity of $\rho_T$ is related
to the above theorem. Note that $\rho_T$ is contractive if
and only if 
$\|f(T)^*\| \leq \|f\|$ by definition of $\rho_T$. Since
$T^* v_i = \bar{z}_i v_i$ we note that $f(T)^*v_i =
\overline{f(z_i)} v_i$, for $1\leq i \leq n$ and $f\in {\mathrm
{Rat}}(\Omega)$. It then follows that
\begin{eqnarray*}
 \|\rho_T(f)^*\|^2 & = & \sup \{\|f(T)^* \big (\sum_{i=1}^n \alpha_i
v_i\big )\|^2 : \alpha_1, \alpha_2, \ldots, \alpha_n \in \C\} \\ 
& = & \sup \{ \|\sum_{i=1}^n \alpha_i \overline{f(z_i)} v_i \|^2 : 
\alpha_1, \alpha_2, \ldots, \alpha_n \in \C\} \\
&=&\sup \{\sum_{i,j=1}^n \alpha_i\bar{\alpha}_j \overline{f(z_i)}f(z_j)
\inner{v_i}{v_j} : \alpha_1, \alpha_2, \ldots, \alpha_n \in \C\}.
\end{eqnarray*}
Therefore, $\|f(T)^*\| \leq \|f\|$ if and only if
$$ \sum_{i,j=1}^n \alpha_i\bar{\alpha}_j \overline{f(z_i)}f(z_j)
\inner{v_i}{v_j} \leq  \sum_{i,j=1}^n
\alpha_i\bar{\alpha}_j\inner{v_i}{v_j}, $$ 
for all $\alpha_1, \alpha_2, \ldots, \alpha_n \in \C$ and all $f \in {\rm Rat}(\Omega)$ with 
$\|f\|\leq 1$.  Thus contractivity of 
$\rho_T$ is equivalent to non-negative definiteness of the matrix
\be \big ( \!\!\big ( (1-\overline{f(z_i)}f(z_j)) K(z_j,z_i) \big
) \!\! \big )_{i,j=0}^n, \ee 
for all  $f \in {\rm Rat}(\Omega)$, $\|f\| \leq 1$. If $\rho_T$ is dilatable then the
theorem above tells us that \be \big ( \!\!\big (
(1-\overline{f(z_i)}f(z_j)) K(z_j,z_i) \big ) \!\! \big
)_{i,j=0}^n  = \big ( \!\!\big (
(1-\overline{f(z_i)}f(z_j))  \langle K_\mathcal{E} (z_j , z_i) x_i , x_j
\rangle \big ) \!\! \big )_{i,j=0}^n. \ee The last matrix is
non-negative definite because $M$ on $H^2_\mathcal{E} (\Omega)$
induces a contractive homomorphism. We therefore see, in this concrete fashion,
that if the homomorphism $\rho_T$ was dilatable then it would be contractive.

The interesting point to note here is that our construction of the
dilation of $\rho_T$ when $T$ is a $2 \times 2$ matrix proves that
the general dilation in that case is of the form $H^2_\bal(\Omega)\otimes\C^2$.

Suppose that the homomorphism $\rho_T$ admits a dilation of the form
\be \label{alphadilation} \rho_{M\otimes I}: \mathcal A(\Omega) \to
\mathcal B(H_\bal^2(\Omega) \otimes\C^n) \ee for some $\bal \in \T^m$, that is, the
multiplication operator $M \otimes I$ on $ H_\bal^2(\Omega)
\otimes\C^n$ is a dilation of $T$. Since the eigenvectors $\{v_1,
v_2, \ldots , v_n\}$ for $T^*$ span $V$ and the set of
eigenvectors of $M^* \otimes I : H^2_\bal(\Omega) \otimes\C^n
\to H^2_\bal(\Omega) \otimes\C^n$ at $z_i$ is the set of vectors
$\{K_\bal(\cdot, z_i)\otimes a_j: a_j \in \C^n, 1\leq j \leq
n\}$ for $1\leq i \leq n$, it follows that any map $\Gamma: V \to
H^2_\bal(\Omega)$ that intertwines $T^*$ and $M^*$ must be
defined by $\Gamma(v_i) = K_\bal(\cdot, z_i)\otimes a_i$ for
some choice of a set of $n$ vectors $a_1, a_2, \ldots ,a_n$ in
$\C^n$.  Now $\Gamma$ is isometric if and only if \be \big (
\!\!\big ( K(z_j ,z_i)\big ) \!\! \big ) =  \big ( \!\!\big (
\inner{v_i}{v_j} \big ) \!\! \big ) = \big ( \!\!\big (
K_\bal(z_j,z_i) \inner{a_i}{a_j}\big ) \!\! \big ). \ee Clearly,
this means that $\big ( \!\!\big ( K(z_j ,z_i)\big ) \!\! \big )$
admits $\big ( \!\!\big ( K_\bal(z_j ,z_i)\big ) \!\! \big )$ as
a factor in the sense that $\big ( \!\!\big ( K(z_j ,z_i)\big )
\!\! \big )$ is the Schur product of $\big ( \!\!\big (
K_\bal(z_j ,z_i)\big ) \!\! \big )$ and a positive definite
matrix, namely, the matrix $A = \big ( \!\!\big (
\inner{a_i}{a_j}\big ) \!\! \big )$.

Conversely, the contractivity assumption on $\rho_T$ does not
necessarily guarantee that $K_\bal$ is a factor of $K$. However, if 
we make this stronger assumption, that is, we assume there exists a positive
definite matrix $A$ such that $\big ( \!\!\big ( K(z_j ,z_i)\big )
\!\! \big ) = \big ( \!\!\big ( K_\bal(z_j,z_i) a_{ij} \big )
\!\! \big )$, where $A = \big ( \!\!\big (a_{ij}\big ) \!\! \big
)$. Since $A$ is positive, it follows that $A = \big ( \!\! \big (
\inner{a_i}{a_j}\big ) \!\! \big )$ for some set of $n$ vectors
$a_1, \ldots , a_n$ in $\C^n$. Therefore if we define the map
$\Gamma: V \to H^2_\bal(\Omega) \otimes \C^n$ to be $\Gamma(v_i)
= K_\bal(\cdot, z_i) \otimes a_i$ for $1\leq i \leq n$ then
$\Gamma$ is clearly unitary and is an intertwiner between $T$ and
$M^*$.  Thus the theorem above has the corollary:

\begin{Corollary} The homomorphism $\rho_T$ is dilatable to a homomorphism
$\tilde{\rho}$ of the form (\ref{alphadilation}) if  the kernel
$K$, as defined in (\ref{k}), is the Schur product of a positive definite matrix $A$ 
and the restriction of $K_\bal$ to the set $\sigma\times \sigma$ for
some $\bal \in \T^m$. 
\end{Corollary}

\noindent
{\sc Acknowledgement: }We thank the referee for pointing out an error in an 
earlier draft of the paper. We also thank him for several comments which helped us in 
our presentation.

\end{document}